\documentclass{article}

\title{Construction of geodesics on Teichm\"uller spaces \\ of Riemann surfaces with $\mathbb Z$ action}
\author{Ryo Matsuda }

\usepackage[ top=30truemm , bottom=30truemm , left=25truemm , right=25truemm]{geometry}
\usepackage{amsmath}
\usepackage{amsthm}
\usepackage{amsfonts}
\usepackage{amssymb}
\usepackage{bm}
\usepackage{color}
\usepackage[dvipdfmx]{graphicx}
\usepackage{framed}
\usepackage{empheq}
\usepackage{enumerate}
\usepackage{fancybox}
\usepackage{fancyhdr}
\usepackage{graphics}
\usepackage{mathrsfs}
\usepackage{ulem}
\usepackage{shadethm}
\usepackage{multicol}
\usepackage[all]{xy}
\usepackage{indentfirst}

\newtheoremstyle{definition}
  {}
  {}
  {\itshape}
  {}
  {\bfseries}
  {\\ \ \ }
  { }
  {}

\theoremstyle{definition}

\newshadetheorem{de}{Definition}

\newtheoremstyle{theorem}
  {}
  {}
  {\itshape}
  {}
  {\bfseries}
  {.\\ \ \ }
  { }
  {}

\theoremstyle{theorem}
\newtheorem{theo}{Theorem}
\newtheorem{copytheo}{Theorem}

\newtheorem{coro}[theo]{Corollary}
\newtheorem{lem}[theo]{Lemma}

\newtheoremstyle{answer}
  {}
  {}
  {}
  {}
  {\itshape}
  {: \\ \ \ }
  { }
  {}
  
\theoremstyle{answer}

\newtheorem*{pr}{Proof}

\newtheoremstyle{ack}
  {}
  {7pt}
  {}
  {}
  {\bfseries}
  {\\ \ \ }
  { }
  {}
  
\theoremstyle{ack}
\newtheorem*{ac}{Acknowlagement.}
\newtheorem{Rem}{Remark}
\newtheorem{exa}{Eg.}


\newenvironment{proofbar}{%
   \MakeFramed {\advance\hsize-\width \FrameRestore}}%
{\endMakeFramed}

{\endMakeFramed}

\newenvironment{lembar}{
   
  \MakeFramed {\advance\hsize-\width \FrameRestore}}
{\endMakeFramed}

\definecolor{lightgray}{rgb}{0.75,0.75,0.75}
\definecolor{problemcolor}{gray}{0.3}
\definecolor{shadecolor}{gray}{0.92}

\allowdisplaybreaks[3] 

\newcommand{\Area}{\mathrm{Area} \ }

\newcommand{\Bel}{\mathrm{Bel}}

\newcommand{\bel}{\mathrm{bel}}

\newcommand{\Cl}{\mathrm{Cl}}

\newcommand{\CZ}{{\mathbb C \setminus \mathbb Z}}

\newcommand{\id} {\operatorname{id}}

\newcommand{\IM}{\mathrm{Im}}

\newcommand{\RE}{\mathrm{Re}}

\newcommand{\Teich}{\mathrm{Teich}}

\newcommand{\dx}{\ dx}
\newcommand{\dtildex}{\ d\tilde x}
\newcommand{\dzeta}{\ d \zeta}

\newcommand{\dz}{\ dz}

\newcommand{\dxdy}{\ dxdy}

\begin{document}

\newcommand{\Addresses}{{
  \bigskip
  \footnotesize

  Ryo Matsuda, \par\nopagebreak
  \textsc{Department of Mathematics, Faculty of Science, Kyoto University, Kyoto 606-8502, Japan}\par\nopagebreak
  \textit{E-mail} : \texttt{matsuda.ryou.82c@st.kyoto-u.ac.jp}
}}

\maketitle


\begin{abstract}
	Teichm\"uller space $\Teich ( R )$ of a Riemann surface $R$ is a deformation space of $R$. 
	In this paper, we prove a sufficient condition for extremality of the Beltrami coefficients 
	when $R$ has the $\mathbb Z$ action.  
	As an application, we discuss the construction of geodesics. 
	Earle-Kra-Krushka\'l proved that 
	the necessary and sufficient conditions for the geodesics connecting $[0]$ and $[\mu]$ to be unique 
	are $\| \mu_0 \|_{\infty} = | \mu_0 | ( z )$ (a.e.$z$)  and ``unique extremality''. 
	As a byproduct of our results, we show that we cannot exclude ``unique extremality''.
	To show the above claim, we  construct a point $[\mu_0]$ in $\mathrm{Teich}(\mathbb C \setminus \mathbb Z)$,  
	satisfying $\| \mu_0 \|_{\infty} = | \mu_0 | ( z )$ (a.e.$z$) and there exists a family of geodesics 
	 $\{ \gamma_\lambda \} _{\lambda \in D}$ connecting $[0]$ and $[\mu_0]$
	with complex analytic parameter, where $D$ is an open set in $l^{\infty}$.	 
\end{abstract}

\section{Introduction}

Teichm\"uller space of a Riemann surface $R$ is the deformation space of complex structure of $R$. 
It is a high-dimensional complex manifold and complete metric space with respect to Teichm\"uller distance, 
and whether its dimension is finite or infinite depends on whether R is of finite or infinite type. 
Since the 1950s, these basic properties have been actively studied, 
as Teichm\"uller space theory with quasiconformal maps was developed by Ahlfors, Bers, and others. 
The detailed history and applications of Teichm\"uller spaces are described in many articles, 
for example, see Ahlfors \cite{QC}, Gardiner \cite{G}, and Hubbard \cite{H}.  
Riemann surfaces of infinite type appear quite naturally, but they are often transcendental. 
As for the geodesics in Teichm\"uller spaces 
treated in this paper, the geodesic connecting two points is unique in the finite-dimensional case, 
but this is generally not the case in the infinite-dimensional case, see section 4 for a detailed explanation. 
Therefore, when dealing with infinite types, it is important to impose a topological or analytical restriction on $R$. 
Hence, in this paper, we deal with the case where $R$ has a $\mathbb Z$-action. 
Imposing conditions on Automorphism of $R$ is, for example, Strebel's problem. 
Ohtake \cite{Oh} showed that in the case of finitely generated Able groups, the lifting by covering the extremal Beltrami coefficients is extremal, see Remark \ref{rem maintheorem1} for detail. 

 
First, we recall the deformation theory of Riemann surfaces. 
Let $R$ be a Riemann surface whose universal covering surface is the upper half plane $\mathbb H$
with a covering map $\pi : \mathbb H \to R$, 
and let be represented by a Fuchsian group $\Gamma$ 
acting on $\mathbb H$ as $R = \mathbb H / \Gamma$. 
Denote by $L^\infty ( \Gamma )$ 
the complex Banach space of the bounded measurable \textit{Beltrami differentials} for $\Gamma$ 
supported on $\mathbb H$. Each element $\mu \in L^\infty ( \Gamma )$ satisfies that 
$\mu \circ \gamma \cdot \overline{ \gamma' } = \mu \cdot \gamma'$ for every $\gamma \in \Gamma$. 
Its open unit ball $\Bel ( \Gamma )$ is the space of Beltrami coefficients for $\Gamma$. 
Each element $\mu \in \Bel ( \Gamma )$, the unique quasiconformal homeomorphism 
$w^\mu : \mathbb H \to \mathbb H$ with Beltrami coefficient $\mu$
 that leaves $0 , 1$ and $\infty$ fixed, is compatible for $\Gamma$. 
Two elements $\mu$ and $\nu$ in $\Bel ( \Gamma )$ are said to be \textit{Teichm\"uller equivalent}
if $w^\mu | _{\mathbb R} \cong w^\nu| _{\mathbb R}$ holds. An element $\mu \in \Bel ( \Gamma )$
is said to be \textit{trivial} if it is Teichm\"uller equivalent to 0. 
The \textit{Teichm\"uller space} $\Teich ( \Gamma )$ of $\Gamma$ 
is the quotient space of $\Bel ( \Gamma )$ by the Teichm\"uller equivalence relation. 
For each $\mu \in \Bel ( \Gamma )$, 
let $[ \mu ]$ denote the point of $\Teich ( \Gamma )$ determined by $\mu$. 
It is known that $\Teich ( \Gamma )$ is a Banach manifold equipped with a complex structure 
and the Teichm\"uller distance. 
The \textit{Teichm\"uller distance} $d_{T ( \Gamma )} ( [ \mu ], [ \nu ] )$
 of two points $[\mu]$ and $[\nu]$ in $\Teich ( \Gamma )$ is defined by 
	\[
		d_{T ( \Gamma )} ( [ \mu ] , [ \nu ] ) 
		: = \inf \tanh^{-1} ( K ( w^{\tilde \nu} \circ (w^{\tilde \mu} )^{-1} ) ), 
	\] 
where $K ( w^{\tilde \nu} \circ (w^{\tilde \mu} )^{-1} )$ represent 
the maximal dilatation, defined in Section 2.1 , of $w^{\tilde \nu} \circ (w^{\tilde \mu} )^{-1}$. 
Let $\Bel(\Gamma)$ and $\Teich(\Gamma)$ be written as $\Bel(R)$ and $\Teich(R)$, respectively. 
See section 2.1 for their detailed relations.

From the above and Section 2.1, 
studying the geometry of Teichm\"uller spaces means clarifying the complex structure of Riemann surfaces. 
In this paper, we work on the former. More specifically, we consider geodesics with respect to the Teichm\"uller distance. 
First of all, the \textit{extremal Beltrami coefficients}, closely related to geodesics, are discussed. 
A Beltrami differential $\mu_0$ in $\Bel ( \Gamma )$ is said to be extremal if 
$\| \mu_0 \|_\infty = \inf \{ \| \mu \|_\infty \mid \mu \in [ \mu_0 ] \}$. 
From the Hamilton--Krushka\'l condition, 
if $\mu_0$ is an extremal Beltrami coefficients, then
the map $[ 0 , 1 ] \ni t \mapsto [ t \mu_0 ]$ is a geodesic from $[ 0 ]$ to $[ \mu_0 ]$, see Section 2.1. 
Hence, we consider extremality under the condition that $R$ has $\mathbb Z$ action. 
The result presented in Theorem \ref{mainre.1} 
is a generalization of Ohtake's results on covering maps and extremality \cite[Theorem 1]{Oh}. 

\begin{copytheo} [Theorem \ref {mainre.1}] 
\begin{oframed}
 
 	Let $R$ be an analytic infinitely Riemann surface that is covering $S$ and the covering group is
	an infinite cyclic group $\langle \gamma \rangle$. 
	Let $\mu$ be in $\Bel ( R )$, 
	which satsfies that there exsist an integrable quadratic holomorphic differntials $\varphi$ with $\int_S |\varphi| = 1$
	and $k \in [ \| \mu \|_{\infty}, 1)$ such that
		\[
			\mu \circ \gamma_n \cdot \frac {\overline { \gamma_n' }} {\gamma_n'} 
			\overset{n \to \infty} {\longrightarrow} 
			k \frac {| \tilde \varphi |} {\tilde \varphi} \ \ \ ( \text {a.e.} \ z \in \omega_0 )
		\]
	where $\tilde \varphi$ is the lift of the element by the covering, $\gamma_n : = \gamma ^{\circ n } $, 
	and $\omega_0$ is fundamental domain for $\langle \gamma \rangle$. 
	Then $\mu$ is extremal on $R$.

\end{oframed}
\end{copytheo}

\begin{Rem} \label{rem maintheorem1}
	 
	Theorem \ref{mainre.1} includes the known fact that
	the lifts of Teichm\"uller Beltrami coefficients on $S$ are extremal in $\Bel ( R )$. 
	There was an important issue involved in these discussions, called Strebel's problem.
	Let $\pi : R \to S$ be a covering of a hyperbolic Riemann surface $S$ 
	with a covering transformation group $\Gamma$. 
	In this situation, there is a natural problem; 
	Under what conditions does $\pi: R \to S$ (or $\Gamma$) satisfy, 
	does the induced map $\pi^\ast : \Bel ( S ) \ni \mu \mapsto \mu \circ \pi \frac {\overline \pi' }{\pi'} \in \Bel ( R )$ 
	preserve extremality? 
	This problem has already been solved. Ohtake \cite[Theorem 1] {Oh} showed that it is sufficient that $\Gamma$
	is a finitely generated Abelian. 
	Finally, McMullen \cite[Theorem 1.1]{Mc} showed that it is equivalent to $\Gamma$ is amenable,
	when $S$ is of finite type. 
	It should be noted,  
	that Ohtake's result can be applied even if $R$ and $S$ are of infinite type, 
	while McMullen's result requires that $S$ be of finite type. 
	Theorem \ref{mainre.1} can be regarded as a generalization of Ohtake's result by weighting. 
	 
\end{Rem}

In the second place, we discuss the construction of a family of geodesics through some two points. 
When $R$ is analytically finite type, 
a geodesic connecting any two points in $\Teich ( R )$ exists uniquely from Teichm\"uller theorem. 
However, when $R$ is of analytically infinite type, it is known to be incorrect. 
Li \cite[Theorem 1]{L1} constructed two points whose geodesics are not unique 
in the universal Teichm\"uller space $\Teich ( \mathbb D )$. 
In general infinite type, 
Tanigawa \cite[Theorem 3.5]{Th} and Li \cite[Theorem 4.1]{L2} 
independently derive sufficient conditions for the non-uniqueness of geodesics. 
In their settings, there exists a family of geodesics through some fixed two points with a complex analytic parameter in 
$\mathbb D$. 
On the other hand, necessary and sufficient conditions for geodesics to be unique are also given in
\cite[Theorem 6]{EKK} and \cite[Theorem 3]{L1}: 
 \begin{enumerate}\renewcommand{\labelenumi}{(\arabic{enumi})}
		\item The Beltrami coefficient $\mu$ is uniquely extremal in $[ \mu ]$, 
		\item $| \mu ( z ) | = \| \mu \|_{\infty}$ \ \ \ a.e. $z \in R$. 
 \end{enumerate}
Condition (1) is called unique extremality, and condition (2) is called absolutely constant. 
Bozin-Lakic-Markovic-Mateljevic \cite[Theorem 10]{BLMM} proved that condition (1) does not imply condition (2). 
Here, 
we prove that condition (2) does not imply condition (1), 
constructing a family of geodesics while maintaining condition (2). 

 \begin{copytheo} [Theorem \ref {mainre. 2}] 
 \begin{oframed}
  
  	If $R$ has $\mathbb Z$ action, 
	then there exists a Beltrami coefficient $\mu$ satisfies the following properties: 
		\vspace{-0.25cm}
		\begin{enumerate}
			\setlength{\parskip}{0cm} 
			\setlength{\itemsep}{0cm} 
				\item $\mu$ is extremal. 
				\item $| \mu |$ is constant
				\item There exists 
				a family of geodesics $\{ g_\lambda \mid \lambda \in \Lambda'  \}$,
				where  $\Lambda' $is an open set in $\mathbb C$, 
				through $[ 0 ]$ and $[ \mu ]$ with a complex analytic parameter.  
				In particular, $\lambda_1 \neq \lambda_2$ implies that 
				the geodesic $[ 0 , 1 ] \ni t \mapsto [ t \mu_{\lambda_1} ]$ 
				is not equal to $ [ 0 , 1 ] \ni t \mapsto [ t \mu_{\lambda_2} ]$.
		\end{enumerate}
		\vspace{-0.4cm}
 \end{oframed}
 \end{copytheo}
This means that 
we cannot exclude the condition (1) ``unique extremality'', which is one of the necessary and sufficient conditions 
for the unique existence of a geodesic given by Earle-Kra-Krushka\'l \cite[Theorem 6]{EKK}.
In addition, if $R = \mathbb C \setminus \mathbb Z$, 
we construct an infinite-dimensional family of geodesics through some fixed two points with a complex analytic parameter in 
an open set of $l^\infty$.

	\begin{copytheo} [Theorem \ref {mainre. 3}] 
	\begin{oframed}
	 
	 	Suppose $R = \mathbb C \setminus \mathbb Z$, then
 		there exists a Beltrami coefficient $\mu$ satisfies following properties: 
		\vspace{-0.25cm}
		 \begin{enumerate}
			\setlength{\parskip}{0cm} 
			\setlength{\itemsep}{0cm} 
				\item $\mu$ is extremal. 
				\item $| \mu |$ is constant. 
				\item There exists  
				a family of geodesics $\{ \lambda_{(\lambda_n )} \mid ( \lambda_n ) \in \Omega  \}$, 
				where  $\Omega$ is an open set in $ l^{\infty} ( \mathbb C )$, 
				through $[ 0 ]$ and $[ \mu ]$ with a complex analytic parameter.  
				In particular, $\lambda_1 \neq \lambda_2$ implies that 
				the geodesic $[ 0 , 1 ] \ni t \mapsto [ t \mu_{\lambda_1} ]$ 
				is not equal to $ [ 0 , 1 ] \ni t \mapsto [ t \mu_{\lambda_2} ]$.
		\end{enumerate}
	 	\vspace{-0.4cm}
	\end{oframed}
	\end{copytheo}

\begin{Rem}
 
 	Historically, the weak version of Theorem B and Theorem C 
	was proved in \cite[Theorem 1]{L1}, \cite[Theorem 4.1]{L2} and \cite[Theorem 3.5]{Th}. 
	However, their results do not consider the condition absolutely constant. Moreover, 
	they constructed only one complex parameter family of geodesics, 
	which is an important extension of the family of geodesics in Theorem C which has
	complex analytic parameters in an infinite-dimensional open set. 
 
\end{Rem}

\begin{ac}
	The author would like to thank Prof. Mitsuhiro Shishikura for 
	giving helpful suggestions, which greatly simplified his previous argument. 
	He also thanks Prof. Hiromi Ohtake for many discussions and helpful advice.
	He is grateful to Yota Maeda for giving him valuable comments on the previous version of this paper. 
	This work was supported by JST, the establishment of 
	University Fellowship Program for the Creation of Innovation in Science and Technology, Grant Number JPMJFS2123.
\end{ac}

\section{Preliminaries}
In this section, we recall some basic facts about quasiconformal maps in the plane, 
integrable holomorphic quadratic differentials, 
extremal Beltrami coefficients, and the general Poincar\'e series. 

\subsection{Quasiconformal maps and Teichm\"uller Theory}

First, we explain 
the quasiconformal maps that play an important role in constructing Teichm\"uller spaces.
Let $K >1$ be a real number. 
A homeomorphic map $f$ that preserves 
the orientation defined by a planar domain $\Omega \subset \mathbb C$
is called a \textit{$K$-quasiconformal} if $f$ has locally integrable distributional derivatives which satisfy 
$| f_{\bar z} | \leq k f_z$ \ a.e. $\Omega$, where $k : = (K -1) / (K + 1)$. 
And $K$ is called \textit{maximal dilatation} of the homeomorphism $f$. 
If $K$ (or $k$) is not specified, the function $f$ is simply called quasiconformal. 
It is well known that $f$ satisfies $f_z \neq 0$ \ a.e. $\Omega$. 
The ratio $\bel ( f ) : = f_{\bar z} / f_z$ is a well-defined measurable function on $\Omega$ with
$\| \bel ( f ) \|_{\infty} \leq k$, and called \textit{complex dilatation} of $f$ or \textit{Beltrami coefficient} of $f$. 
Conversely, from the measurable Riemann mapping theorem, for every bounded measurable function $\mu$ on $\mathbb C$
 with $\| \mu \|_\infty <1$, there exists a unique quasiconformal map $f^\mu$ with Beltrami coefficient $\mu$
 that leaves $0, 1$ and $\infty$ fixed, see \cite[Chapter 5]{QC}. 
 Hereafter, a quasiconformal map between Riemann surfaces means 
 that each composite of local coordinates is a quasiconformal map and the dilatation of the map $\mu$ satisfies
 $\| \mu \|_{\infty} < 1$. 
 
Next, we briefly  explain that Teichm\"uller space $\Teich(\Gamma )$ can also be regarded as a deformation space 
of the complex structure of Riemann surface $R$, using quasiconformal maps in the plane. 
Let $L^\infty (R)$ denote the complex Banach space of 
all bounded measurable $(-1, 1 )$-forms on $R$, where $\mu$ is a bounded measurable $(-1, 1 )$-form
means $\mu = \mu dz / d\bar z$. 
The elements contained in its open unit ball $\Bel(R)$ are called the Beltrami coefficients on $R$. 
Note that an element $\mu$ of $L^\infty (R)$ lifts to an element $\hat \mu$ of $L^\infty (\Gamma)$ by the formula 
$\hat \mu = \mu \circ \pi \cdot \frac { \bar{\pi'} } {\pi'}$. The map $\mu \mapsto \hat \mu$ is an isometry between $L^\infty (R)$
to $L^\infty (\Gamma)$. 
Let $\mu$ be an element in $\Bel(R)$, then $\hat \mu$ is an element in $\Bel (\Gamma)$. 
Since $w^{\hat \mu}$ is a homeomorphism and compatible with $\Gamma$, 
$\Gamma^\mu$ is a discrete group acting on $\mathbb H$, 
hence $R^\mu : = \mathbb H / \Gamma^\mu$ becomes a Riemann surface again. 
Moreover $w^{\hat \mu}$ induces a quasiconformal map $f^\mu$ between $R \to R^\mu$. 
From classical results, for two elements $\mu$ and $\nu$ in $\Bel(R)$, 
$\hat \mu$ is Teichm\"uller equivalent to $\hat \nu$ if and only if there exists a conformal map $c : R^\mu \to R^\nu$
such that $( f^{\nu} )^{-1} \circ c \circ f^\mu$ is homotopic to $\id_R$ relative ideal boundary of $R$. 
The Teichm\"uller space of $R$ is $\Bel(R)$ factored by Teichm\"uller equivalence relation, and denoted by $\Teich(R)$. 
 $\Teich (R)$ and $\Teich (\Gamma)$ are naturally isomorphic, see \cite[Proposition 1 in \S 5]{G}.
 
\subsection{Quadratic differentials and Beltrami differntials}

Let $Q(R)$ denote the complex Banach space of all integrable holomorphic quadratic differentials 
for $\Gamma$. The space $Q(\Gamma)$ is identified with the cotangent space of $\Teich(\Gamma)$
at $[ 0 ]$. The complex dimension of $\Teich(\Gamma)$ and $Q(\Gamma)$ are finite, if and only if
$\Gamma$ is finitely generated of the first kind. In particular, The dimension is equal to $3g - 3 + n$, 
where $g$ is the genus and $n$ is a number of punctures of $R$. 

A Beltrami differential $\mu_0$ in $\Bel ( \Gamma )$ is said to be extremal if 
$\| \mu_0 \|_\infty = \inf \{ \| \mu \|_\infty \mid \mu \in [ \mu_0 ] \}$. 
It is known that a Beltrami differential $\mu$ is extremal if and only if there exists a sequence
$( \varphi_j )_{j \in \mathbb N}$ of elements 
in $Q ( \Gamma )$ with unit integrable norm such that 
$\| \mu \|_{\infty} = \lim_{n \to \infty} \RE \int_{\mathbb H / \Gamma} \mu \varphi_n$
, it is called the Hamilton--Krushkal condition, see \cite[Theorem 1 and Theorem 6 in \S 6.7]{G}. 
Such a sequence is called a Hamilton sequence for $\mu$. 
From the Hamilton--Krushkal condition, if $\mu$ is extremal, 
we see that the following map is geodesic from $[ 0 ]$ to $[ \mu ]$ with respect to the Teichm\"uller distance: 
	\[
		[ 0 , 1 ] \ni t \mapsto [ t \mu ]. 
	\]
A Beltrami differential of the form $z | \varphi | / \varphi $ 
with some $z \in \mathbb D$ and
$\varphi \in Q ( \Gamma )$ is called Teichm\"uller Beltrami differential, 
where $\mathbb D$ is the unit disk in the complex plane.
It is obvious that Teichm\"uller Beltrami differentials are extremal. 
It is known that a point $[ \mu ]$ in $\Teich ( \Gamma )$ which contains a Teichm\"uller Beltrami differential is
the unique extremal differential in equivalence class $[ \mu ] $, see \cite[Theorem 2 in \S 6.2, Theorem 3 in \S6.3]{G}. 

Let $S$ be a Reimann surface with the universal covering transformation group $\Gamma'$ acting $\mathbb H$. 
Suppose that $R$ is a covering surface of $S$, then $\Gamma$ is a subgroup of $\Gamma'$. 
For any element $\varphi$ in $Q ( R )$, we set
	\[
		\Theta ( \varphi ) ( z ) : = \sum_{g \in \Gamma / \Gamma'} \varphi \circ g \cdot ( g ' ) ^ 2(z), 
	\]
where $\Gamma / \Gamma'$ is the coset. 
It is said to be a (general) \textit{Poincar\'e series}. The series is absolutely and uniformly convergence on $\omega_0$, 
where $\omega_0$ is a fundamental domain in $\mathbb H$ for $\Gamma'$. 
Moreover $\Theta ( \varphi )$ is an element in $Q ( S )$. 
The linear operator $\Theta: Q ( R ) \to Q ( S )$ is bounded and surjective.

\setcounter{theo}{0}
\section{Extremality under $\mathbb Z$ action}

We discuss the extremality of the Beltrami coefficient, 
which is closely related to the geodesics of Teichm\"uller space. 
The Hamilton-Krushka\'l condition is a well-known necessary and sufficient condition: 
A Beltrami coefficient $\mu$ on a Riemann surface $R$ is extremal if and only if 
		\[
			\| \mu \|_{\infty} = \sup \left\{ \int_\omega \mu \varphi 
			\middle | \varphi \in Q(R), \| \varphi \|_{Q(R)} = 1 \right\}, 
		\]
where $\omega$ is a fundamental domain. 
The proof is founded in many articles, for example, see\cite[Theorem 1  and Theorem 6 in \S 6.7]{G}.
Now we prove the following which plays a critical role in the discussion below.  

\begin{theo} [ ] \label{mainre.1}
\begin{oframed}
	
	Let $R$ be an analytic infinitely Riemann surface that is covering $S$ and the covering group is
	an infinite cyclic group $\langle \gamma \rangle$. 
	Let $\mu$ be in $\Bel ( R )$, 
	which satisfies that there exsist $\varphi \in Q ( S )$ with $\| \varphi \|_{Q ( S )} = 1$
	and $k \in [ \| \mu \|_{\infty}, 1)$ such that
		\[
			\mu \circ \gamma_n \cdot \frac {\overline { \gamma_n' }} {\gamma_n'} 
			\overset{n \to \infty} {\longrightarrow} 
			k \frac {| \tilde \varphi |} {\tilde \varphi} \ \ \ ( \text {a.e.} \ z \in \omega_0 )
		\]
	where $\tilde \varphi$ is the lift of the element $\varphi$ in $Q ( S )$, $\gamma_n : = \gamma ^{\circ n } $, 
	and $\omega_0$ is fundamental domain for $\langle \gamma \rangle$. 
	Then $\mu$ is extremal. 
	
\end{oframed}
\end{theo}

\begin{pr}
\begin{proofbar}

	Let $\omega_0$ is a fundamental domain for $\langle \gamma \rangle$ in $R$, 
	$\omega_n : = \gamma_n ( \omega_0 )$, and $D_n : = \bigcup_{0 \leq j \leq n} \omega_j$.
	Since the general Poincar\'e series $\Theta : Q (R) \to Q (S)$ is surjective, there is an element $f$ in $Q ( R )$ 
	such that $\tilde \varphi = \Theta ( f )$ on $\omega_0$. 
	Set 
		\[
			F_n : = 
			\frac{1} {\|  \Theta ( f ) \|_{Q (S)} } 
			\sum_{-n \leq j \leq 0} f \circ \gamma_j \cdot ( \gamma_j' ) ^2, \ \ \  
			\varphi_n : = \frac {1} {\| F_n \|_{Q ( R )}} F_n.
		\]
	We will prove that $( \varphi_n )$ is a Hamilton sequence for $\mu$. Namely, we prove that 
		\[
			\left| \int_R \mu \varphi_n \right| \overset{n \to \infty} {\longrightarrow} k. 
		\]
	Since the following inequality: 
			\begin{eqnarray*}
				k & \geq & \left| \int_R \mu \varphi_n \right| \\
				& = & 
				\frac {1} {\| F_n \|_{Q(R)}} \left| \int_{D_n} \mu \tilde \varphi 
				+ \int_{D_n} \mu ( F - \tilde \varphi )
				+ \int_{R \setminus D_n} \mu F_n \right| \\
				& \geq & 
				\frac {1} {\| F_n \|_{Q(R)}}
				\left( \left| \int_{D_n} \mu \tilde \varphi \right| 
				- \left| \int_{D_n} \mu ( F - \tilde \varphi ) \right| - 
				\left| \int_{R \setminus D_n} \mu F_n \right| 
				\right) \\
				& = & 
				\frac {n + 1} {\| F_n \|_{Q(R)}} 
				\left\{
				\frac{1} {n + 1}  \left|\int_{D_n} \mu \tilde \varphi \right| - 
				\frac {1} {n + 1} \left( \left| \int_{D_n} \mu ( F - \tilde \varphi ) \right| +
				\left| \int_{R \setminus D_n} \mu F_n \right|  \right)
				\right\}, 
			\end{eqnarray*}
	it is sufficient to show that the following four sequences converge as shown below:  
		\[
			\lim_{n \to \infty} \frac {1} {n + 1}\int_{D_n} \mu \tilde \varphi = k, \ \ 
			\lim_{n \to \infty} \frac {1} {n + 1} \int_{D_n} | \tilde \varphi - F_n |  = 0, 
		\]
		\[
			\lim_{n \to \infty} \frac {1} {n + 1}\left| \int_{R \setminus D_n} \mu F_n \right| = 0, \ \ 
			\lim_{n \to \infty} \frac {\| F_n \|_{Q(R)}} {n + 1} = 1. 
		\]

	First, let us $ x_n : = \int_{\omega_n} \mu \tilde \varphi$, then we get 
		\begin{eqnarray*}
    			x_n
    			& = & 
    			\int_{\omega_0} 
    			\mu \circ \gamma_n \cdot \tilde \varphi \circ \gamma_n \cdot | \gamma_n' | ^2
    			 =  \int_{\omega_0} 
    			\mu \circ \gamma_n \cdot \tilde \varphi \cdot \frac{1} {(\gamma_n ')^2} \cdot | \gamma_n' | ^2 
    			 =  \int_{\omega_0} 
    			\mu \circ \gamma_n \cdot \frac {\overline { \gamma_n' }} {\gamma_n'} \cdot \tilde \varphi \\
    			& \overset{n \to \infty} {\longrightarrow} & 
    			\int_{\omega_0}k \frac {| \tilde \varphi |} {\tilde \varphi} \tilde \varphi. 
    			= k \| \varphi \|_{Q(S )} = k, 
    		\end{eqnarray*}
	because $\mu$ is Beltrami differential and $\tilde \varphi$ is holomorphic quadratic differntial.  
	Hence, from Ces\'aro mean, it follows that	
		\[
			\frac {1} {n + 1} \int_{D_n} \mu \tilde \varphi = \frac {1} {n + 1} \sum_{0 \leq j \leq n} x_j
			\overset{n \to \infty} {\longrightarrow} k. 
		\]

	Second, by the definition of $\tilde \varphi$ and $F_n$, we calculate
		\begin{eqnarray*}
				\frac {1} {n + 1} \int_{D_n} | \tilde \varphi - F_n | 
				& \leq &
				\frac {1} {n + 1} 
				\sum_{0 \leq l \leq n} 
				\int_{\omega_l} \sum_{j > n \ \text {or} \ j < 0 } | f \circ \gamma_j| \cdot | ( \gamma_j' )^2 | \\
				& = & 
				\frac {1} {n + 1} 
				\sum_{0 \leq l \leq n} 
				\int_{\omega_0} \sum_{j > n \ \text {or} \ j < 0 } | f \circ \gamma_j \circ \gamma_l |
				\cdot | ( \gamma_j' )^2 | | ( \gamma_l' )^2 | \\
				& = & \frac {1} {n + 1} 
				\sum_{0 \leq l \leq n} 
				\int_{\omega_0} \sum_{m - l > n \ \text {or} \ m - l < 0 } | f \circ \gamma_m |
				\cdot | ( \gamma_m' )^2 | \overset{n \to \infty} {\longrightarrow}  0. 
		\end{eqnarray*}

	Third, in the same way as above, the following calculation holds: 
		\begin{eqnarray*}
				\frac {1} {n + 1} 
				\int_{R \setminus D_n} \sum_{-n \leq j \leq 0} | f \circ \gamma_j \cdot ( \gamma_j' ) ^2 | 
				& = & 
				\frac {1} {n + 1} 
				\sum_{0 > l \ \text {or} \ l > n} 
				\int_{\omega_l} \sum_{-n \leq j \leq 0} | f \circ \gamma_j \cdot ( \gamma_j' ) ^2 | \\
				& = & 
				\frac {1} {n + 1} 
				\sum_{0 > l \ \text {or} \ l > n} 
				\int_{\omega_l} \sum_{-n \leq j -l \leq 0} | f \circ \gamma_j \cdot ( \gamma_j' ) ^2 | 
				\overset{n \to \infty} {\longrightarrow} 0
		\end{eqnarray*}

	Finally, we will prove that 
		\[
			\left|  \| \Theta(f) \|_{Q ( S )} - \frac {1} {n + 1} \int_{D_n} | F_n | \right| 	
			\overset{n \to \infty} {\longrightarrow} 0, 
		\]
	it implies that 
		\[
			\frac {1} {n + 1} \int_{D_n} | F_n | 
			\overset{n \to \infty} {\longrightarrow} \| \Theta(f) \|_{Q ( S )}. 
		\]
	From the definition of Poincar\'e series, it is seen that 
		\begin{eqnarray*}
				\left|  \| \Theta(f) \|_{Q ( S )} - \frac {1} {n + 1} \int_{D_n} | F_n | \right| 	
				& = &
				\frac {1} {n + 1} 
				\left|  
				( n + 1 ) \| \Theta(f) \|_{Q ( S )} 
				- \int_{D_n} \left| \sum_{-n \leq j \leq 0} f \circ \gamma_j \cdot ( \gamma_j' ) ^2 \right| \right| \\ 			
				& = & 
				\frac {1} {n + 1} 
				\left|  
				\sum_{0 \leq l \leq n } \int_{\omega_0} | \Theta ( f ) | 
				- \sum_{0 \leq l \leq n } \int_{\omega_0} 
				\left| \sum_{-n \leq j - l \leq 0} f \circ \gamma_j \cdot ( \gamma_j' ) ^2 \right| \right| \\
				& \leq & 
				\frac {1} {n + 1} 
				\sum_{0 \leq l \leq n } 
				\int_{\omega_0}
				\left| \Theta ( f ) -
				\sum_{-n \leq j - l \leq 0} f \circ \gamma_j \cdot ( \gamma_j' ) ^2 \right| \\
				& \leq & 
				\frac {1} {n + 1} 
				\sum_{0 \leq l \leq n } 
				\int_{\omega_0}
				\sum_{-n > j - l \text \ {or} \ j -l >0}\left|  f \circ \gamma_j \cdot ( \gamma_j' ) ^2 \right| 
				\overset{n \to \infty} {\longrightarrow} 0
			\end{eqnarray*}
	Thus, we get
		\[
			\frac {\| \Theta ( f ) \|_{Q(S)}\| F_n \|_{Q(R)}} {n + 1} 
				= 
				\frac {1} {n+1} \left( \int_{D_n} + \int_{R \setminus D_n} \right) | F_n | 
				\overset{n \to \infty} {\longrightarrow} \| \Theta ( f ) \|_{Q(S)} + 0. 
		\]
	Summing up the above calculations, we can show that the sequence $(\varphi_n)$ in $Q(S)$ satisfies
		\[
			\| \varphi_n \|_{Q(S)} = 1 \ \ ( \forall n ), \ \ 
			\lim_{n \to \infty} \left| \int_R \mu \varphi_n \right| = k = \|\mu\|_\infty.
		\]
	Hence, It can be seen that $\mu$ is extremal from the Hamilton-Krushka\'l condition. 
	\qed
	
\end{proofbar}
\end{pr}

We compose an example that satisfies the assumption of Theorem \ref{mainre.1}. 

\begin{exa}
	
	Let us look at the infinitely analytic Riemann Surface $R : = \mathbb C \setminus \mathbb Z$, 
	then $R$ has action by infinite cyclic group generated by $z \mapsto z + 1$. 
	Set $S : = R / \langle \gamma : z \mapsto z + 3 \rangle$, and 
	$\omega_0 : = \{ z \in R \mid 0 < \RE z \leq 3 \}$.  
	Since $S$ is conformal equivalent to 
	five punctured sphere, it follows that $\dim_{\mathbb C} Q ( S )$ is equal to $2$ 
	from Riemann -Roch Theorem. 
	Let $\varphi^r$ and $\varphi^l$ from $Q ( S )$ be linear independent, 
	and let $(a_n)$ and $( b_n)$ be sequences in $[ 0 , 1 )$ which satisfies that $a_n + b_n \neq 0$ and
		\[
			\lim_{n \to \infty} a_n = 1, \lim_{n \to - \infty} a_n = 0, 
			\lim_{n \to \infty} b_n = 0, \lim_{n \to - \infty} b_n = 1. 
		\]
	We define Beltrami coefficient $\mu$ on $R$ as follows: 
		\[
			\mu ( = \mu_{(a_n, b_n)} ) : = 
			k \frac {| a_n \varphi^r + b_n \varphi^l |} {a_n \varphi^r + b_n \varphi^l}
			 \ \ \ ( z \in \omega_n), 
		\]	
	where $k$ is a constant in $[ 0 ,1 )$, $\omega_n : = \gamma^{ \circ n} ( \omega_0 )$, 
	and we consider $\varphi^r$ and $\varphi^l$ as elements in $Q ( S )$,identifying $\omega_j$ with $S$.
	A simple computation shows that $\mu$ satisfies the assumption of Theorem \ref{mainre.1}. 
	Indeed, 
		\[
			\mu \circ \gamma_n \cdot \frac {\overline { \gamma_n' }} {\gamma_n'} ( z )
			= 
			k \frac {| a_n \varphi^r + b_n \varphi^l |} {a_n \varphi^r + b_n \varphi^l} ( z ) 
			\overset{n \to \infty} {\longrightarrow} k \frac {|\varphi^r |} {\varphi^r} ( z )
			\ \ \ ( \text{for all} \ z \in \omega_0 ). 
		\]

	For instance, if we take $a_n = 1, b_n = 0$, then $\mu_{(a_n, b_n)}$ can be regarded as the lift of $\varphi^r$. 
	Therefore Theorem \ref{mainre.1} is a generalization of Ohtake's result(\cite[Theorem 1]{Oh}). 

\end{exa}

\section{Construction of geodesics}

The uniqueness of geodesics plays an important role in considering the geometry of Teichm\"uller space. 
In the finite-dimensional case, the geodesic connecting any two points is unique. 
By contrast, in the infinite-dimensional case, there exists a point $[ \mu ]$, which satisfies that 
the geodesic connecting $[ 0 ]$ and $[ \mu ]$ is not unique. 
Of course, there are two points whose geodesics are unique. 

The first result on this issue is that Li \cite[Theorem 1]{L1} constructed an example in $\Teich ( \mathbb D )$ . 
In general infinite-dimensional Teichm\"uller spaces, 
Li \cite[Theorem 4.1]{L2} and Tanigawa \cite[Theorem 3.5]{Th} independently constructed examples. 
On the other hand, Li \cite[Theorem 3]{L1} and Earle-Kra-Krushka\'l \cite[Theorem 6]{EKK} proved that 
the necessary and sufficient condition for a geodesic connecting $[ 0 ]$ and $[ \mu ]$ 
to be unique is that the following two conditions hold: 
	\begin{enumerate}\renewcommand{\labelenumi}{(\arabic{enumi})}
		\item The Beltrami coefficient $\mu$ is uniquely extremal in $[ \mu ]$, 
		\item $| \mu ( z ) | = \| \mu \|_{\infty}$ \ \ \ a.e. $z \in R$. 
	\end{enumerate}
In \cite{L1}, Li asked if condition (1) implies condition (2). However,
for this problem, Bozin-Lakic-Markovic-Mateljevic gave a counterexample \cite[Theorem 10]{BLMM}. 
That is to say, they constructed that 
$\mu$ is uniquely extremal and $\mu$ is not absolutely constant in a certain Riemann Surface. 
From the above, we consider how many geodesics can be constructed while preserving conditions (2). 
Namely, 
we configure a family of Beltrami coefficients $\{ \mu_\lambda \}_{\lambda \in \Lambda}$ which satisfies
the following properties: 
	\vspace{-0.1cm}
	\begin{enumerate}
	\setlength{\parskip}{0cm} 
	\setlength{\itemsep}{0.2cm} 
	\renewcommand{\labelenumi}{\alph{enumi}).}
		\item For all $\lambda \in \Lambda$, $\mu_\lambda$ is extremal. 
		\item For all $\lambda_1, \lambda_2 \in \Lambda$, $\mu_{\lambda_1}$. 
		is Teichm\"uller equivalent to $\mu_{\lambda_2}$. 
		\item $\lambda_1 \neq \lambda_2$ implies $[ 0 , 1 ] \ni t \mapsto [ t \mu_{\lambda_1} ]$
		and $[ 0 , 1 ] \ni t \mapsto [ t \mu_{\lambda_2} ]$
		are distinct.  
		\item There exists $\lambda_0 \in \Lambda$ 
		such that $| \mu_{\lambda_0} |( z )= \| \mu_{\lambda_0} \|_{\infty}$ \ \ \ (a.e. $z \in R$).
	\end{enumerate}
 The fact that the tangent vectors at the origin of each geodesic 
 are different is a sufficient condition for condition c). Such a characterization is given by Li \cite{L2}. 
 
 \begin{theo}  [{\cite[Theorem 3.1]{L2}}]
 \begin{oframed}
 	
	Let $\mu_1$ and $\mu_2$ be extremal Beltrami differentials in $[ \mu ] \in \Teich ( \Gamma )$. 
	A sufficient condition for two geodesics 
	$[ 0 , 1 ] \ni t \mapsto [ t\mu_1]$ and $[ 0 , 1 ] \ni t \mapsto [ t\mu_2]$  
	to be different is that there exists $\varphi \in Q ( \Gamma )$ such that 
		\[
			\int_{\mathbb H / \Gamma } ( \mu_1 - \mu_2 ) \varphi \neq 0. 
		\]	
\vspace{-0.4cm}
 \end{oframed}
 \end{theo}
 
To achieve our goal, we fix one open subset U of R and construct it in the following form: 
	\[
		\mu_\lambda : = \chi_{U} \tau_\lambda + \chi_{R \setminus U} \mu,  
	\]
 where $\chi_U$ is the characteristic function of $U$ and $\{ \tau_\lambda \}_{\lambda \in \Lambda}$ 
 is a family of Beltrami coefficients on $U$. 
 In particular, by choosing $\mu$ well, conditions a) and b) are satisfied,
 and by constructing  $\{ \tau_\lambda \}_{\lambda \in \Lambda}$ well, conditions c) and d) are satisfied. 
 
 \subsection{Construction of the family $\{ \tau_\lambda \}$}
 
 First, to define $\{ \mu_\lambda \}$, we construct a family of Beltrami coefficients on $U$. 
 
  \begin{lem} \label{ trivial qc on strip }
 \begin{lembar}
 	
	Let $T_{R} : = \{ \zeta \in \mathbb C \mid  0 < \IM \zeta < R \}$, 
	and $\Lambda : = \{ \lambda \in \mathbb C \mid | \IM \lambda | < 1, \RE \lambda > 0 \}$ \ ($R \in ( 0 , \infty )$). 
	There exists a family of Beltrami coefficients $\{ \widehat{ \tau_\lambda } \}_{\lambda \in \Lambda}$ on $T_{R}$
	which satisfies the following properties: 

		\vspace{-0.1cm}
		\begin{itemize}
		\setlength{\parskip}{0cm} 
		\setlength{\itemsep}{0.2cm} 
			\item For all $\lambda \in \Lambda$, 
			$\widehat{ \tau_\lambda }$ is Teichm\"uller equivalent to $0$ in $\Bel ( T_{R})$. 
			\item $\lambda \in \mathbb R \cap \Lambda$ implies that $| \widehat{ \tau_\lambda } |$ is constant. 
			\item $\lambda_1 \neq \lambda_2$ implies that 
			$\int_{ \mathbb T_{R} / \langle \zeta \mapsto \zeta + 2 \pi \rangle } 
			( \widehat{ \tau_{\lambda_1} } - \widehat{ \tau_{\lambda_2} } ) 1 \dxdy \neq 0$. 
			\item 
			$\widehat{ \tau_{\ast} } : \Lambda \ni \lambda
			\mapsto \widehat{ \tau_\lambda } \in L_\infty (T_{R} )$ is holomorphic map.  
		\end{itemize}
 \end{lembar}
 \end{lem}
 
 \begin{pr}
 \begin{proofbar}
  
  	Consider the following a self-affine map $F_\lambda : T_{R} \to T(r)$, 
		\[
			F_\lambda : T_{R} ( \xi + i\eta \mapsto ) : =
			\begin{cases}
				\xi + \lambda\eta + i\eta & ( 0 < \eta \leq \frac {R} {2}), \\
				\xi + \lambda ( R - \eta ) + i\eta & (\frac {R} {2} \leq \eta < R).
			\end{cases}
		\]
	Let $\widehat{ \tau_\lambda }$ be a beltrami coefficient of $F_\lambda$. 
	Since $F_\lambda |_{\partial T_{R}} = \id_{\partial T_{R}}$, 
	$\widehat{ \tau_\lambda }$ is Teichm\"uller equivalent to $0$ in $\Bel ( T_R )$. 
	The explicit calculation of $\tau$ is as follows: 
		\begin{eqnarray}
			\widehat{ \tau_\lambda } = 
			\begin{dcases}
				\frac {i\lambda} {2 - i\lambda} & \left(0 < \eta \leq \frac {R} {2}\right) , \\
				\frac {-i\lambda} {2 + i\lambda} & \left(\frac {R} {2} \leq \eta < R\right) .
			\end{dcases}
		\end{eqnarray}
	Thus, the remaining properties are obvious. \qed
		
 \end{proofbar}
 \end{pr}

 Since $\pi : T_R \ni \zeta \mapsto e^{i\zeta} \in \mathbb A_{r_0} : = \{ z \in \mathbb C \mid r_0 < | z | < 1 \}$, 
 where $r_0 : = e^{-R}$, is a covering map, and $F_c$ satisfies that $F_\lambda ( \zeta + 2 \pi ) = F_\lambda ( \zeta ) + 2 \pi$, 
 $F_\lambda$ induces a self quasiconfoermal map $f_c$ of $\mathbb A_{r_0} $. 
 Let $\nu_c$ be a Beltrami coefficient of $f_c$, then the pullback of $\nu_c$ by $\pi$ coincide with $\widehat{ \tau_\lambda }$. 
 Note that the pullback of $\dz^2 / z^2$ by $\pi$ is $\dzeta^2$ and 
 that the third complement is the coupling of the Beltrami coefficients and the holomorphic quadratic differentials 
 the following proposition holds: 
 
 \begin{coro} \label{ trivial qc on annulli }
 \begin{lembar}
  
 	Set $\mathbb A_{r} : = \{ z \in \mathbb C \mid r < | z | < 1 \}$, where $ r \in [ 0 , 1 )$. 
	There exists a family of Beltrami coefficients $\{ \nu_\lambda \}_{\lambda \in \Lambda}$ on $\mathbb A_{r} $
	which satisfies the following properties: 
		\vspace{-0.1cm}
		\begin{itemize}
		\setlength{\parskip}{0cm} 
		\setlength{\itemsep}{0.2cm} 
			\item For all $\lambda \in \Lambda$, 
			$\nu_\lambda$ is Teichm\"uller equivalent to $0$ in $\Bel ( \mathbb A_{r} )$. 
			\item $\lambda \in \mathbb R \cap \Lambda$ implies that $| \nu_\lambda |$ is constant. 
			\item $\lambda_1 \neq \lambda_2$ implies that 
			there exists holomorphic function $g$ on $\mathbb A_{r}$ such that
			$\int_{\mathbb A_{r}} ( \nu_{\lambda_1} - \nu_{\lambda_2} ) g \dxdy \neq 0$. 
			\item $\nu_{\ast} : \Lambda \ni \lambda \mapsto \nu_\lambda \in L_\infty ( \mathbb A_{r} )$ 
			is a holomorphic map.  
		\end{itemize}

 \end{lembar}
 \end{coro}

 \begin{pr}
 \begin{proofbar}
 
 	The first claim is more obvious from the properties of $F_\lambda$. 
	Explicitly calculating $\nu_\lambda$ is as follows: 
		\begin{eqnarray} \label {eq : anu tri belcoe 1}
			\nu_\lambda = 
			\begin{dcases}
				\frac {-i\lambda} {2 + i\lambda } \frac {z} {\bar z} & (r < | z | \leq \sqrt {r}), \\
				\frac {i\lambda} {2 - i\lambda }\frac {z} {\bar z} & (\sqrt {r} \leq | z | < 1).
			\end{dcases}
		\end{eqnarray}
	A simple computation shows that $\lambda \in \mathbb R \cap \Lambda$ implies $| \tau_\lambda |$ is constant. 
	The remaining properties follow from the previous discussion. \qed
	
%
 
 \end{proofbar}
 \end{pr}
 
 \begin{Rem}
 
 	Since the lifts of the holomorphic quadratic differentials can distinguish $\nu_\lambda$ in the infinitesimal sense is 
	$\dzeta^2$, if the coefficient of $z^{-2}$ of the holomorphic function $h$ on $\mathbb A_{r}$ is zero, then 
		\[
			\int_{\mathbb A_{r_0}} \nu_c h \dxdy = 0. 
		\]
	For later discussion, we consider the coupling between quadratic differentials and Beltrami coefficients on $T_R$. 
	Let $\tilde \varphi_\ast$ be a holomorphic quadratic differential on $T_R$, then $\tilde \varphi_\ast$ satisfies 
		\[
			\tilde \varphi_\ast ( \zeta ) = \tilde \varphi_\ast \circ \iota ( \zeta ) \cdot ( \iota' ( \zeta ) ) ^2
			= \tilde \varphi_\ast ( \zeta + 2 \pi ), 
		\]
	where $\iota ( \zeta ) = \zeta + 2 \pi$. Namely, $\tilde \varphi_\ast$ is regarded as a periodic function. 
	In $\mathbb T_R : = T_R / \langle \iota \rangle$, Fourier expansion yields: 
		 \[
			\int_{\mathbb T_R} \nu_\lambda \tilde \varphi_\ast 
			= \tilde \varphi_\ast (0) \Area ( \mathbb T_R )  
			\left( \frac {i\lambda} {2 - i c} - \frac {i\lambda} {2 + i \lambda} \right)
			=  \tilde \varphi_\ast (0) \Area ( \mathbb T_R ) \frac {\lambda^2} {4 + \lambda^2}. 
		\]
	In the following, we will use the case $h(z) = \frac {1} {z - \alpha}$, where $| \alpha | \leq r_0$, and we calculate it.
	Since  the coefficient of $z^{-2}$ of $h$ is 
		\[
			\frac {1} {2 \pi i} \int_{| z | = r} \frac {z} {z - \alpha} \dz
			= \frac {1} {2 \pi i} \left( \int_{| z | = r} 1 + \frac {\alpha} {z - \alpha} \dz \right) = \alpha, 
		\]
	we get 
		\[
			\iint_{\mathbb A_{r}} \nu_\lambda \frac {1} {z - \alpha} 
			= - \frac {2 \pi \alpha} {i} \Area(\mathbb A_{r}) \frac {-\lambda^2} {4 + \lambda^2}.
		\]
  \end{Rem}
 
 Taking the limit of $r_0$ to $0$, we can construct a family of Beltrami coefficients on $\mathbb D$ 
 which satisfies the same properties. 
 In the following, we describe different construction methods using elliptic functions. 
 
 \begin{theo} [ ] \label { trivial qc on disk }
 \begin{oframed}
 	
	Let $\Lambda \subset \mathbb C$ take the same as above. 
	There exists a family of Beltrami coefficients $\{ \tau_\lambda \}_{\lambda \in \Lambda}$ on $\mathbb D$
	which satisfies the following properties: 
		\vspace{-0.1cm}
		\begin{enumerate}
		\setlength{\parskip}{0cm} 
		\setlength{\itemsep}{0.2cm} 
		\renewcommand{\labelenumi}{\alph{enumi}').}
			\item For all $\lambda \in \Lambda$, 
			$\tau_\lambda$ is Teichm\"uller equivalent to $0$ in $\Bel ( \mathbb D )$. 
			\item $\lambda \in \mathbb R \cap \Lambda$ implies that $| \tau_\lambda |$ is constant. 
			\item $\lambda_1 \neq \lambda_2$ implies that 
			there exists holomorphic function $g$ on $\mathbb D$ such that
			$\int_{\mathbb D} ( \tau_{\lambda_1} - \tau_{\lambda_2} ) g \dxdy \neq 0$. 
			\item $\tau_{\ast} : C \ni c \mapsto \tau_\lambda \in L^\infty ( \mathbb D )$ is holomorphic map.  
		\end{enumerate}
		\vspace{-0.4cm}

 \end{oframed}
 \end{theo}
 
 \begin{pr}
 \begin{proofbar}
 
 	Since Riemann mapping theorem, for all $R \in ( 0, \infty )$, there uniquely exists 
	${t} \in ( 0 , 1 )$ 
	which satisfies that there exists a covering map $\rho : T_R \to D_{t}$, 
	where $D_{t}$ is the set $\mathbb D \setminus [ - t, t]$. 
  	In particular, we take $\tilde \rho$ such that $\rho ( \{ \zeta | \IM \zeta = R \} ) = [ -t, t] $, where
	$\tilde \rho$ is extension of $\rho$ to the closure of $T_R$ to the closure of $D_{t}$.
	
	The $\rho$ is the elliptic function theta satisfies the following differential equation: 
		\[
			(\rho') ^2 
			= (\rho - t ) ( \rho + t ) \left( \rho - \frac {1} {t} \right)\left( \rho + \frac {1} {t} \right).
		\]

	We construct the family of Beltrami coefficients on $D_{t}$
	from the  family of Beltrami coefficients $\{ F_c \}$ constructed by the Lemma \ref { trivial qc on strip }
	using $\rho$. In fact, $F_\lambda$ is compatible with $\zeta \mapsto \zeta + 2 \pi$, 
	it induces a self quasiconformal map of $D_t$:  
		\[
			\xymatrix{
			T_R \ar[d]^-{\rho} \ar[r]^-{F_\lambda}  &T_R \ar[d]^-{\rho} \\
			D_{t} \ar[r]^-{\exists \tilde {f_\lambda}} & D_{t} \ar@{}[lu]|{\circlearrowright}
			}
		\]
	Let $\tau_\lambda$ be a beltrami coefficient of $\tilde { f_\lambda }$. 
	It is easy to verify that the family of Beltrami coeffinietns $\{ \tau_\lambda \}_{\lambda \in \Lambda}$ 
	satisfies a'), b'), and d' ). We claim that the family satisfies  c'). 
	
	Recall the remark of Corollary \ref { trivial qc on annulli }, 
	we prove that 
	there exists a holomorphic function $\varphi_\ast$ 
	such that the pullback $\tilde \varphi_\ast : = \rho^\ast ( \varphi_\ast \dz^2 )$ satisfies that 
	$\tilde \varphi_\ast ( 0 ) \neq 0$. 
	From simple calculations, we obtain 
		\begin{eqnarray*}
			\tilde \varphi_\ast ( 0 )
			& = & 
			\int_{- \pi} ^\pi \tilde \varphi_\ast ( \tilde x ) \dtildex 
			= 2 \int_{0} ^\pi \tilde \varphi_\ast ( \tilde x ) \dtildex \\
			& = &
			2 \int_{0}^\pi \varphi_\ast \circ \rho \cdot (\rho')^2 \dtildex \\
			& = &
			2 \int_{0} ^{t} \varphi_\ast ( x ) \cdot (\rho') \dx \\
			& = &
			2 \int_{0}^{t} 
			\varphi_\ast (x) \cdot \sqrt{(x - t) ( x + t) 
			\left( x - \frac {1} {t} \right)\left( x + \frac {1} {t} \right)} \dx. 
		\end{eqnarray*}
	Note that integrand is positive on Integral interval, $\tilde \varphi_\ast ( 0 ) \neq 0$ 
	under the condition which $\varphi_\ast( 0 ) \neq 0$ and $R$ is sufficiently large. 
	\qed

 \end{proofbar}
 \end{pr}

 \begin{Rem}
  
  	In \cite{K}, Reich constructed a family of Beltrami coefficients 
	$\{ \tau_\lambda \}_{\lambda \in C}$
	on $\mathbb D$ which 
	satisfies the following properties, where $C : = \{ \lambda \in \mathbb C \mid | \lambda | < 1/2 \}$: 
		\begin{itemize}
		\setlength{\parskip}{0cm} 
		\setlength{\itemsep}{0cm} 
			\item For all $\lambda \in C$, 
			$\tau_\lambda$ is Teichm\"uller equivalent to $0$ in $\Bel ( \mathbb D )$. 
			\item $\lambda_1 \neq \lambda_2$ implies that 
			there exists holomorphic function $g$ on $\mathbb D$ such that
			$\int_{\mathbb D} ( \tau_{\lambda_1} - \tau_{\lambda_2} ) g \dxdy \neq 0$. 
			\item $\tau_{\ast} : C \ni \lambda \mapsto \tau_\lambda \in L^\infty ( \mathbb D )$ is holomorphic map.  
		\end{itemize}
	However, his example does not have elements whose absolute values are constant. 
 \end{Rem}

 \subsection{For sums of Beltrami coefficients whose supports are disjoint}
 
Next, let us consider condition b) explained at the beginning of this section. 
Let $R_1$ be an infinitely analytic subsurface of $R$ whose boundary is the union of relative compact analytic curves in $R$,
and $R_2 : = R \setminus R_1$. 
In general, $R_2$ is probably not connected. Then, we consider the following function: 
	\[
		\Bel ( R_1 ) \ni \tau \mapsto \tau + \mu \in \Bel (R), 
	\]
where $\mu$ is Beltrami coefficitent on $R_2$ and 
$\mu + \tau$ is extended to be identically $\mu$ in $R_2$. 

\begin{theo} [ ] \label {lem : trivial + mu is trivial}
\begin{oframed}
 
 	The above map induces the holomorphic map: 
		\[
			\Teich( R_1 ) \ni [\tau] \mapsto [ \tau + \mu ] \in \Teich (R). 
		\]
	In other words, if $\tau_1$ is Teichm\"uller equivalent to $\tau_2$ on $R_1$, 
	$\mu + \tau_1$ is Teichm\"uller equivalent to $\mu + \tau_2$ on $R$. 
 
\end{oframed}
\end{theo}

\begin{pr}
\begin{proofbar}
	
	First, we prove that $ \tau_1 + 0$ is Teichm\"uller equivalent to $\tau_2 + 0$ on $R$. 	
	Because $\tau_1$ is Teichm\"uller equivalent to $\tau_2$ on $R_1$, there exists a conformal map
	$c : f^{\tau_1} ( R_1 ) \to f^{\tau_2} ( R_1 )$ such that $( f^{\tau_2} ) ^ {-1} \circ c \circ f^{\tau_1}$ is 
	homotopic to $\id_{R_1}$ in the following sence: There exists 
	a homotopy $( g_t : R_1 \to R_1 )_{t \in [ 0 , 1 ]}$ which extend continuously to the border of $R_1$ 
	such that $g_0 = \id_R$, $g_1 = ( f^{\tau_2} ) ^ {-1} \circ c \circ f^{\tau_1}$ and $g_t |_{\partial R_1} = \id$. 
		\[
			\xymatrix{
				f^{\tau_1} ( R_1 ) \ar[dd]^-c & \quad & \quad & f^{0+\tau_1} ( R ) \ar[dd]^-{\hat c} \\
				\quad & R_1 \ar[lu]^-{f^{\tau_1}} \ar[ld]^-{f^{\tau_2}} \ar@{^{(}-_>}[r]
				& R  \ar[ru]^-{f^{0+\tau_1}} \ar[rd]^-{f^{0 +\tau_2}} & \quad \\
				f^{\tau_2} ( R_1 ) & \quad & \quad & f^{0+\tau_2} ( R ) 
			}
		\]
	We consider the following maps: 
		\[
			\hat {c} ( p ) : = 
			\begin{cases}
				f^{\tau_2 + 0} \circ ( f^{\tau_2} ) ^{-1} \circ c \circ f^{\tau_1} \circ ( f^{\tau_1 + 0} ) ^{-1} ( p )
				& p \in f^{\tau_1 + 0} ( R_1 ) \\
				f^{0 + \tau_2} ( f^{\tau_1 + 0} ) ^{-1} ( p )
				& p \in f^{\tau_1 + 0} ( R_2 )
			\end{cases},
			\ \ \ 
			\hat g_t ( p ) : = 
			\begin{cases}
				g_t ( p ) & p \in R_1 \\
				\id & \in R_2
			\end{cases}. 
		\]	
	Note that $g_t |_{\partial R_1} = \id$ implies that $\hat {g} _t$ is self-continuous map of $R$.

	$\hat c$ is conformal and 
	$( f^{\tau_2 + 0} ) ^ {-1} \circ \hat c \circ f^{\tau_1 + 0}$ is homotopic to $\id_R$. Indeed, 
	$\bel ( f^{\tau_1} \circ ( f^{ \tau_1 + 0} ) ^{-1} ) = 0$, since support of $\tau_1$ is contained in $R_1$.  
	Moreover $g_t$ is the homotopy which join $\id_{R_1}$ and $( f^{\tau_2} ) ^ {-1} \circ c \circ f^{\tau_1}$, thus
	$\hat g_t$ join  $\id_{R}$ and $( f^{ \tau_2 + 0} ) ^ {-1} \circ \hat c \circ f^{\tau_1 + 0}$. 
	
	Next  we prove that $\tau_1 + \mu$ is Teichm\"uller equivalent to $\tau_2 + \mu$ on $R$. 	
		\[
			\xymatrix{
				f^{\tau_1 + 0} ( R ) \ar[dd]^-{ \hat c } 
				& \quad & \quad & \quad & f^{\mu+\tau_1} ( R ) \ar[dd]^-{\hat { \hat c }} \\
				\quad & R \ar[lu]^-{f^{\tau_1}} \ar[ld]^-{f^{0 + \tau_2}} \ar[rr]^{f^{\mu + 0}}
				\ar[rrru]^-{f^{\mu+\tau_1}} \ar[rrrd]_{f^{\mu +\tau_2}}
				& \quad & f^{\mu + 0} (R)  \ar[ru]_{h^{\tau_1}} \ar[rd]^-{h^{\tau_2}} & \quad \\
				f^{0 + \tau_2} ( R) & \quad & \quad & \quad & f^{\mu+\tau_2} ( R ) 
			}, 
		\]
 	where  $h^{\tau_j} : = f^{\mu + \tau_j} \circ ( f^{\mu + 0} )^{-1}$ ( $j = 1,2 )$ and 
		\begin{eqnarray*}
			\hat {\hat c} &: = & 
			f ^{\mu + \tau_2} \circ ( f^{0 + \tau_2} ) ^{-1} \circ \hat c \circ f^{0 + \tau_1} \circ ( f^{\mu + \tau_1} ) ^{-1}.
		\end{eqnarray*}
	Similarly to the above proof, 
	we can show that $\mu + \tau_1$ and $\mu + \tau_2$ are Teichm\"uller equivalence in $f^{\mu + 0} ( R )$. 
	In detail, $	\hat { \hat g } _t : = ( f ^{\mu+0} ) ^{-1} \circ g_t \circ f^{\mu+0}$ is the homotopy which joins 
	$\id_{f^{\mu + 0} ( R )}$ and $( h^{\tau_2} ) ^{-1} \circ \hat {\hat c} \circ h^{\tau_1}$. 
	Therefore $\mu + \tau_1$ is Teichm\"uller equivalent to $\mu + \tau_2$ on $R$. \qed

\end{proofbar}
\end{pr}

\begin{Rem}

Taniguchi \cite{Tm} and Maitani \cite{Mf} 
considered a subsurface $R'$ of the Riemann surface $R$ such that
if two Beltrami coefficients on the subsurface $R'$ which are Teichm\"uller equivalent, 
the homotopy appearing in the Teichmuller equivalence on $R'$ can be extended to $R \setminus R'$ by the identity map. 
Using their results, Tanigawa proved that 
 if $R_1$ is simply connected, the induced map is injective,
 see \cite[Lemma3.3, Lemma 3.4]{Th} for detail. 

\end{Rem}

\subsection{Main results}

 Before starting the main results, 
 let us prepare a lemma on the construction of the integrable holomorphic quadratic differentials, 
 which is necessary for the proof of them. 
 
 \begin{lem}\label {lem : quad diff pole}
 \begin{lembar}
  
  	There exists an integrable holomorphic quadratic differential on $R$ which has a pole of order $1$ at a puncture in $R$. 

 \end{lembar}
 \end{lem}
 
 \begin{pr}
 \begin{proofbar}
  
  	Set $R' : = R \cup \{ a \}$. Let $R'$ be represented by a Fuchsian group $\Gamma'$ 
	whose a covering map $\pi : \mathbb D \to R'$ satisfies $\pi ( 0 ) = a$.  
	Note that $\psi (z) : = 1 / z$ is an integrable meromorphic function on $\mathbb D$. 
	Hence, $\Psi : = \Theta ( \psi )$ has properties of the claim. \qed
 \end{proofbar}
 \end{pr}
 
 Through the above discussion, we obtain the main theorem in this paper. 
 
 \begin{theo} [ ] 
 \begin{oframed}
  
  	If $R$ has action of infinite cyclic group $\langle \gamma \rangle$, then 
	there exists a Beltrami coefficient $\mu$ satisfies the following properties: 
		\vspace{-0.25cm}
		\begin{enumerate}
			\setlength{\parskip}{0cm} 
			\setlength{\itemsep}{0cm} 
				\item $\mu$ is extremal.
				\item $| \mu |$ is constant.
				\item There exists a family of Beltrami coefficients
				$\{ \mu_\lambda \}_{\lambda \in \Lambda}$ such that 
					\begin{enumerate}
						\item $\mu \in \{ \mu_\lambda \}_{\lambda \in \Lambda}$,
						\item all included Beltrami coefficients are Teichm\"uller equivalence,
						\item $\Lambda \ni \lambda \mapsto \mu_\lambda \in \Bel ( R )$ is holomorphic map,
						\item $\lambda_1\neq \lambda_2$ implies 
						$\mu_{\lambda_1}$ is not infinitesimal equivalent to $\mu_{\lambda_2}$, and
						\item there exists a domain $U \subset R$
						such that $\mu_c |_{R \setminus U} = \mu|_{R \setminus U}$ and 
						$\Cl (U) \cap \{ a_n \} \neq \emptyset$. 
					\end{enumerate}
			\end{enumerate}
		In particular, 
		$U$ is conformal equivalent to $\mathbb D$, $\mathbb D \setminus \{ 0 \}$, or $\mathbb A$.

 \end{oframed}
 \end{theo}
 
\begin{pr}
\begin{proofbar}
 
 	Compose $\tilde \mu$ in the same way as in Eg 1. 
	Namely, let $\varphi^r$ and $\varphi^l$ from $Q ( R / \langle \gamma \rangle )$ be linear indipendent, 
	and let $X_n : = (a_j, b_j )$ be sequences in $\mathbb R_{\geq 0}^2 \setminus \{ ( 0, 0 ) \} $ with
	$\lim_{n \to \infty} X_n = (0 , 1)$. Using these, we define $\psi$ to be $a_n \varphi^r + b_n \varphi^l$, and
		\[
			\tilde \mu (z) : = t_0 \frac {| \psi_n |} {\psi_n} ( z ) \ \ \ ( z \in \omega_n ), 
		\]   
	where $t_0$ is a constant in $[ 0 ,1 )$, $\omega_j : = \gamma^{ \circ j} (\omega_0)$, 
	and we consider $\varphi^r$ and $\varphi^l$ as elements in $Q ( R / \langle \gamma \rangle ) )$, 
	identifying $\omega_j$ with $R / \langle \gamma \rangle )$. \\

	\noindent \fbox{ \noindent (I).  Configuration where $U$ is conformal equivalent to $\mathbb D$} 
	
	Let $l$ be an analytic simple closed curve 
	such that One of the connected components that $l$ separates conformal equivalent to $\mathbb D$. 
	Then we denote by $U$ the connected component and by $h$ a Riemann map from $\mathbb D$ to $U$. 
	Using the family of Beltrami coefficients constructed in Theorem \ref { trivial qc on disk }, 
	we construct a family of trivial Beltrami coefficients
	$\{ \mu_\lambda : = \widetilde{\tau_\lambda} + \chi_{R \setminus U} \widetilde { \mu } \}_{\lambda \in \Lambda}$, 
	where $\tilde { \tau }: = h_\ast (\mu) $ is the push-forward of the Beltrami coefficients on $\mathbb D$ by $h$, 
	considering $\tilde { \tau } = 0$ on $R \setminus U$. 
	In the following, we show that satisfies $\{ \mu_\lambda \}$
	the claim in the theorem. 
	
	Set $\mu : = \widetilde{\tau_{\lambda_0}} + \chi_{R \setminus U} \widetilde { \mu }$, where $\lambda_0$ satisfies that
		\[
			\left| \frac {-i \lambda_0 / 2} { 1 +i \lambda_0 / 2 } \right| = t_0.
		\]
	It is easy to see that $| \mu |$ is constant and $\mu$ is extremal. 
	Thus, (1) and (2) are shown.
	
	The conditiions b) and d) are trivial, so we will prove a) and c).  
	For all $\lambda \in \Lambda$, by Theorem \ref {lem : trivial + mu is trivial}, 
	$\mu_\lambda$ is Teichm\"uller equivalent to $ \chi_{R \setminus U} \widetilde { \mu }$. 
	Moreover, suppose that $\lambda_1 \neq \lambda_2$. Let  $\varphi$ be in $Q ( R )$  satisfies  
	$\varphi ( h^{-1} ( 0 ) ) \neq 0$. Then we get 
		\[
			\int_R ( \mu_{ \lambda_1 } - \mu_{\lambda_2} ) \varphi = 
			\int_{U} ( \widetilde { \tau_{\lambda_1} } - \widetilde { \tau_{\lambda_2} } ) \varphi \neq 0.  
		\]
	Hence, a) and c) are proved.
	
	\noindent \fbox{ \noindent (II). Configuration where $U$ is conformal equivalent to $\mathbb D^\ast$}
	
	Let $U$ be a topological punctured disk in $R$ with the analytic boundary. We denote by $h$ a conformal map 
	from $U$ to $\mathbb D ^\ast$ and $\alpha$ a puncture contained by $U$. 
	Replace $\{ \tau \}$ in the above configuration (I) with $\{ \nu \}$ 
	constructed in Corollary \ref { trivial qc on annul } applying $r = 0$ to form a family. That is, 
	we consider 
	$\{ \mu_\lambda : = \widetilde{\nu_\lambda} + \chi_{R \setminus U} \widetilde { \mu } \}_{\lambda \in \Lambda}$. 
	Since the only nontrivial condition is condition c, we show this. 
	
	Let $\lambda_1 \neq \lambda_2 $. Since Lemma \ref{lem : quad diff pole}, 
	there exists $\varphi \in Q ( R )$ which has a poke of order 1 at $\alpha$. Therefore 
		\[
			\int_R ( \mu_{\lambda_1} - \mu_{\lambda_2 )} \varphi
			= \int_{U} ( \widetilde{\nu_{\lambda_1}} -  \widetilde{\nu_{\lambda_1}} ) \varphi \neq 0.
		\]

	\noindent \fbox{ \noindent (III). Configuration where $U$ is conformal equivalent to $\mathbb A$} 
	
	Use the same domain $U$, the Riemann map $h$, and the puncture $\alpha$ as the above configuration (II). 
	Next, let $D$ be an annulus contained in $U$ such that $\partial D$ contains the puncture $\alpha$. 
	Replace $\{ \nu \}$ in the above configuration with $\{ \tau \}$ 
	constructed in Corollary \ref { trivial qc on annul } to form a family. That is, 
	we consider 
	$\{ \mu_\lambda : = \widetilde{\nu_\lambda} + \chi_{R \setminus U} \widetilde { \mu } \}_{\lambda \in \Lambda}$. 
	The only nontrivial condition c) is shown in the same way as the above configuration (II). 
	The above observation completes the proof. 
	\qed

\end{proofbar}
\end{pr}
 
 From this theorem, the following assertions about geodesics are derived. 
 
 \begin{theo} [ ] \label{mainre. 2}
 \begin{oframed}
  
  	If $R$ has $\mathbb Z$ action, 
	then there exists a Beltrami coefficient $\mu$ that satisfies following properties: 
		\vspace{-0.25cm}
		\begin{enumerate}
			\setlength{\parskip}{0cm} 
			\setlength{\itemsep}{0cm} 
				\item $\mu$ is extremal. 
				\item $| \mu |$ is constant.
				\item There exists 
				a family of geodesics $\{ g_\lambda \mid \lambda \in \Lambda'  \}$,
				where  $\Lambda' $ is an open set in $\Lambda$, 
				through $[ 0 ]$ and $[ \mu ]$ with a complex analytic parameter.  
				In particular, $\lambda_1 \neq \lambda_2$ implies that 
				the geodesic $[ 0 , 1 ] \ni t \mapsto [ t \mu_{\lambda_1} ]$ 
				is not equal to $ [ 0 , 1 ] \ni t \mapsto [ t \mu_{\lambda_2} ]$.
		\end{enumerate}
  		\vspace{-0.4cm}
 \end{oframed}
 \end{theo}
 
 \begin{pr}
 \begin{proofbar}
  
  	To prove the theorem, we only need to construct 
	a family consisting of extreme Beltrami coefficients that satisfies the conditions.
	Consider the family of Beltrami coefficients $\{ \mu_\lambda \}_{\lambda \in \Lambda}$ and 
	the Beltrami coefficient $\mu$ constructed in the previous theorem.
	Put 
		\[
			\Lambda' : = \left\{ \lambda \in \Lambda \middle | 
			\left| \frac {-i \lambda_0 / 2} { 1 + i \lambda_0 / 2 } \right| < t_0 ( = \| \mu \| ) \right\}. 
		\]
	Then, for all $\lambda \in \Lambda'$, 
	$\| \mu_\lambda \| $ is not greater than $t_0$. Hence, $\mu_\lambda$ is extremal from Theorem \ref {mainre.1}. 
	Finally, we set $g_\lambda : [ 0 , 1 ] \ni t \mapsto [ t \mu_\lambda ]$. \qed
  
 \end{proofbar}
 \end{pr}
 
 As an above theorem, we can prove the following claim:
 
 \begin{coro}
 \begin{lembar}
  
  	We cannot exclude ``unique extremality'', which 
	is one of the necessary and sufficient conditions for the geodesics connecting $[0]$ and $[\mu]$ to be unique
	in \cite[Theorem 6]{EKK}. 
  
 \end{lembar}
 \end{coro}
 
 Since the technical difficulty is the composition of integrable holomorphic quadratic differentials which distinguish geodesics, 
 more geodesics can be constructed if $R$ is a subdomain of the plane, 
 using what is known about integrable holomorphic quadratic differentials in detail.
 For example, more geodesics can be constructed when considered to $\mathbb C \setminus \mathbb Z$. 
 	
\begin{theo} [ ] \label{mainre. 3}
\begin{oframed}
 	
	Suppose $R = \mathbb C \setminus \mathbb Z$, then
 	there exists a Beltrami coefficient $\mu$ satisfies following properties: 
		 \vspace{-0.25cm}
		 \begin{enumerate}
			\setlength{\parskip}{0cm} 
			\setlength{\itemsep}{0cm} 
				\item $\mu$ is extremal. 
				\item $| \mu |$ is constant. 
				\item There exists  
				a family of geodesics $\{ \lambda_{(\lambda_n)} \mid (\lambda_n) \in \Omega  \}$, 
				where  $\Omega$is an open set in $ l^{\infty} (\mathbb C)$, 
				through $[ 0 ]$ and $[ \mu ]$ with a complex analytic parameter.  
				In particular, $\lambda_1 \neq \lambda_2$ implies that 
				the geodesic $[ 0 , 1 ] \ni t \mapsto [ t \mu_{\lambda_1} ]$ 
				is not equal to $ [ 0 , 1 ] \ni t \mapsto [ t \mu_{\lambda_2} ]$.
		\end{enumerate}
		 \vspace{-0.4cm}
\end{oframed}
\end{theo}
 
\begin{pr}
\begin{proofbar}
 
 	Set $S : = R / \langle z \mapsto z + 3 \rangle$, and make $\mu_0$ one of the components of Eg1.
	Let $K > 0$ satisfies
		\[
			\left | \frac { - \frac {iK} {2}  } { 1 - \frac {iK} {2}  } \right| = \| \mu_0 \|.
		\] 
	For this $K$, let denote by $\Omega : = \{ (\lambda_j ) \in l^\infty \mid \lambda_j \in \Lambda, | \lambda_j | \in [ 0, K] \}$. 
	
	First, for any $j \in \mathbb N$, let denote by 
		\[
			\mathbb A_j : = \Delta \left ( 3j + \frac {3} {2} ; \frac{1} {2} + \frac {1} {2^{|j|}} \right) 
			\setminus \Delta \left( 3j + \frac {3} {2} ; \frac {1} {2} \right), \ \ \ 
			U : =   \bigcup_{j \in \mathbb Z} \mathbb A_j.
		\]
	Moreover, denote  $\alpha_j = 3j + 1, \beta_j : = 3j + 2$ , and $U  \bigcup_{j \in \mathbb Z} \mathbb A_j$. 
	
	Next, for each $( \lambda_n ) \in \Omega$, we define 
		\[
			\mu_{(\lambda_j)} : = 
			\chi_{(\CZ) \setminus U} \mu_0 + \sum_{j \in \mathbb Z} \nu_{\lambda_j} \chi_{\mathbb A_j}, 
		\]
	useing $\nu_j$ is constructed in Corollary \ref { trivial qc on annulli }, applying $\mathbb A_j$. 
	Note that each satisfies the following equality: 
		\[
			\iint_{\mathbb A_j} \nu_j \frac{1} {z - \alpha_j} 
			= - \frac{2 \pi \alpha_j} {i} \Area ( \mathbb A_{j}) \frac{A_j^2} {1- A_j^2}, \ \ \
			\iint_{\mathbb A_j} \nu_j \frac{1} {z - \beta_j} 
			=-\frac{2 \pi \beta_j} {i} \Area ( \mathbb A_{j}) \frac{A_j^2} {1 - A_j^2}, 
		\]
	where $A_j : = \frac {- i \lambda_j} {2}$. Moreover if $\varphi$ is a holomorphic on 
	some simply connected domain containing the closure of $\Delta ( 3j + 3 / 2 ; 1 )$, then 
	$\int_{\mathbb A_j} \nu_j \varphi = 0$. 
	
	We show that the constructed 
	$\{ \mu_{(\lambda_j)} \}$ are all Teichm\"uller equivalences. 
	Let $\bf{\nu} : = \sum_{j \in \mathbb Z} \nu_{\lambda_j} \chi_{\mathbb A_j}$, then 
	$f^{\nu}$ is homotopic to $\id_{\CZ}$ relative $\mathbb Z$, therefore 
	 $f^{\chi_{(\CZ) \setminus U} \mu_0} \circ f^{\bf {\nu} }$ is homotopic to 
	 $f^{\chi_{(\CZ) \setminus U} \mu_0} $. Since the intersection of supports of
	 $\chi_{(\CZ) \setminus U} \mu_0$ and $\bf {\nu}$ is empty, we get
		\[
			\bel ( f^{\chi_{(\CZ) \setminus U} \mu_0} \circ f^{\bf {\nu} } ) 
			= \chi_{(\CZ) \setminus U} \mu_0 + \sum_{j \in \mathbb Z} \nu_{\lambda_j} \chi_{\mathbb A_j}.
		\]
	That is, $\{ \mu_{(\lambda_j)} \}$ is Teichm\"uller equivalent to  $\chi_{(\CZ) \setminus U} \mu_0$. 
	
	Finally, 
	we see that $[ 0 , 1 ] \ni t \mapsto [ t \mu_{(\lambda_j)} ]$ gives different geodesics, 
	using contra position. In detail,  we prove that if $[ 0 , 1 ] \ni t \mapsto [ t \mu_{(\lambda_j)} ]$ is equal to 
	$[ 0 , 1 ] \ni t \mapsto [ t \mu_{(\tilde \lambda_j)} ]$, then $\lambda_j = \tilde \lambda_j$ 
	for all $j \in \mathbb N$. Given $L > 0$. 
	Set 
		\[
			\varphi_{\alpha_j} : = \frac{1} { (z - \alpha_j)(z - 3L) ( z - 6L)},
			\varphi_{\beta_j} : = \frac{1} { (z - \beta_j)(z - 3L) ( z - 6L)}. 
		\]
	Let $( \cdot_j )_{| j | \leq L}$ be a finite string of length $2 L + 1$ consisting of $\alpha$ and $\beta$, 
	and  $\varphi_{( \cdot )} : = \sum_{| j | \leq L} \varphi_{\cdot_j}$. 
	Since $\mu_{(\lambda_j)}$ and $\mu_{(\tilde \lambda_j})$ give the same geodesic, applying 
	the finite string $( \alpha, \alpha, \cdots, \alpha)$, we get 
		\begin{eqnarray*}
			0 & = & \int_{\CZ} ( \mu_{(c_j)} - \mu_{( \tilde c_j )} ) \varphi_{(\cdot)} \\
			& = & 
			- \frac{2 \pi } {i} \sum_{-L \leq j \leq L} \alpha_j
			\Area ( \mathbb A_{j}) \left( \frac{A_j^2} {1- A_j^2}
			-\frac{\tilde A_j^2} {1- \tilde A_j^2} \right)
			+ \sum_{| j | > L} \int_{\mathbb A_j} ( \mu_{(c_j)} - \mu_{( \tilde c_j )} ) \varphi_{(\cdot)}. 
		\end{eqnarray*}
	Besides, for any $-L \leq J \leq L$ , we replace 
	the $J$-th of the string $( \alpha, \alpha, \cdots, \alpha)$ 
	from $\alpha$ to $\beta$. By drawing the equation above against the equation obtained 
	by applying the holomorphic function corresponding to the substituted string to the different geodesics, 
	we obtain 
		\[
			\frac{A_J^2} {1- A_J^2} -\frac{\tilde A_J^2} {1- \tilde A_J^2} = 0.
		\]
	Hence, $\lambda_J = \tilde \lambda_J $. \qed
	 
\end{proofbar}
\end{pr}

\Addresses

\end{document}